\documentclass[12pt]{amsart}
\usepackage{times,amsmath,amsthm,amsfonts,amssymb}
\usepackage{mathrsfs}
\usepackage{fullpage}

\title[Interpolation and Sampling in higher dimensions]{Interpolation and
Sampling Hypersurfaces for the Bargmann-Fock space in higher dimensions}

\author{Joaquim Ortega-Cerd\`a}
\address{Departament de Matem\` atica Applicada i An\` alisi
\newline \indent
Universitat de Barcelona, Gran Via 585, 08071-Barcelona, Spain}
\email{jortega@ub.edu}

\author{Alexander Schuster}
\address{Department of Mathematics \newline\indent
San Francisco State University, San Francisco, CA 94132}
\email{schuster@sfsu.edu}

\author{Dror Varolin}
\address{Department of Mathematics\newline\indent
University of Illinois at Urbana-Champaign, Urbana, IL 61801}
\email{dror@math.uiuc.edu}

\date{}

\newcommand{\BF}{\mathfrak{BF}}

\newcommand{\ms}{\medskip}

\newcommand{\rmk}{\noi {\bf Remark.}\ \ }
\newcommand{\noi}{\noindent}

\newcommand{\ifix}{\!\!\!\!}
\newcommand{\iifix}{\!\!\!\!\!\!\!\!}

\newcommand{\cl}{{\mathcal L}}
\newcommand{\co}{{\mathscr O}}
\newcommand{\so}{{\mathscr O}}

\newcommand{\ii}{\sqrt{-1}}

\newcommand{\C}{{\mathbb C}}
\newcommand{\D}{{\mathbb D}}

\newcommand{\di}{\partial}
\newcommand{\dbar}{\bar \partial}
\newcommand{\ddbar}{\partial\bar \partial}

\newcommand{\vp}{\varphi}
\newcommand{\ve}{\varepsilon}

\renewcommand{\Re}{\operatorname{Re}}
\newcommand{\vol}{\operatorname{vol}}
\newcommand{\dist}{\operatorname{dist}}
\newcommand{\area}{\operatorname{Area}}

\newcommand{\lf}{\mathfrak{bf}}

\newtheorem{thm}{Theorem}
\newtheorem{mthm}[thm]{Theorem}
\newtheorem{lem}{Lemma}[section]
\newtheorem{prop}[lem]{Proposition}
\newtheorem{cor}[lem]{Corollary}
\newtheorem*{claim}{Claim}

\newtheorem*{defn}{Definition}

\theoremstyle{remark}

\begin{document}
\begin{abstract}
We study those smooth complex hypersurfaces $W$ in $\C ^n$ having the property
that all holomorphic functions of finite weighted $L^p$ norm on $W$ extend to
entire functions with finite weighted $L^p$ norm.  Such hypersurfaces are
called interpolation hypersurfaces.
We also examine the dual problem of finding all sampling hypersurfaces, i.e.,
smooth hypersurfaces $W$ in $\C ^n$ such that any entire function with finite
weighted $L^p$ norm is stably determined by its restriction to $W$.

We provide sufficient geometric conditions on the hypersurface to be
an interpolation and sampling hypersurface.  The geometric conditions
that imply the extension property and the restriction property are
given in terms of some directional densities.
\end{abstract}

\maketitle

\section*{Introduction}
\addtocounter{section}1

Let $\omega = \ii \di \dbar |z|^2$ denote the standard Euclidean
form in $\C ^n$. Fix a smooth closed complex hypersurface $W \subset \C ^n$
and a plurisubharmonic function $\vp$ such that for some contants
$C,C'>0$,
$$
C\omega \le \ii \di \dbar \vp \le C' \omega
$$
in the sense of currents.  For brevity, such an estimate will sometimes be
denoted $\ii \di \dbar \vp \simeq \omega$.  

For $p\in[1,\infty)$, let
\[
\BF_{\varphi}^p(\C ^n) :=\left \{ F \in \co (\C^n)\ ;\ \int _{\C
^n}|F|^p e^{-p\varphi} \omega ^n < +\infty \right \}
\]
and
\[
\lf _{\varphi}^p (W) := \left \{ f \in \co (W)\ ;\ \int _W |f|^p
e^{-p\varphi} \omega ^{n-1} <+\infty \right \}
\]
denote the generalized Bargmann-Fock spaces of weighted $L^p$
holomorphic functions on $\C ^n$ and $W$ respectively.
When $p=+\infty$ we replace the integrals by suprema.  The classical
Bargmann-Fock space corresponds to the case $\vp (z) = |z|^2$.

\begin{defn}
Let $W$ be a uniformly flat smooth hypersurface in $\C ^n$.
(See section~\ref{unif-flat-section}.)
\begin{enumerate}
\item  We say $W$ is an interpolation hypersurface if for each
$f\in \lf^p _{\varphi}  (W)$ there exists $F \in \BF^p
_{\varphi} (\C ^n)$ such that $F |W = f$.
\item  We say $W$ is a sampling hypersurface if there is a
constant $M=M(p,W)$ such that for all
$F \in \BF^p_{\varphi}
(\C ^n)$,
\begin{eqnarray}\label{samp-ineq}
\frac{1}{M} \int _{\C ^n} |F|^p e^{-p\varphi} \omega ^{n} \le \int
_W |F|^p e^{-p\varphi} \omega ^{n-1} \le M  \int _{\C ^n} |F|^p
e^{-p\varphi} \omega ^{n}
\end{eqnarray}
when $p < +\infty$, or a similar estimate involving suprema in
place of integrals when $p=\infty$.
\end{enumerate}
\end{defn}

The goal of this paper is to find geometric sufficient conditions
for a uniformly flat hypersurface $W$ to be interpolating or
sampling.  A key concept is given in the following definition.

\begin{defn}Let $T \in \co (\C ^n)$ be a holomorphic function such that $W =
T^{-1} (0)$ and $dT$ is nowhere zero on $W$.  For any $z\in\C^n$ and
any $r>0$ consider the (1,1)-form
\[
\Upsilon _W(z,r) := \sum _{i ,\bar j =1} ^n \left ( \frac{1}
{\vol(B(z,r))}\int _{B(z,r)} \frac{\di ^2\log |T|}{\di \zeta ^i
\di \bar \zeta ^j} \omega ^n(\zeta) \right ) \ii dz^i \wedge d\bar
z ^j.
\]
\end{defn}

\medskip

\rmk {\it Clearly the definition of $\Upsilon _W (z,r)$ is
independent of the choice of the function $T$ defining $W$.
Moreover, if $\Theta_W$ is the current of integration associated to
$W$ then $\Upsilon _W$ is the average of $\Theta _W$ in a ball of
center $z$ and radius $r>0$:
\[
\Upsilon _W = \Theta_W * \frac {\mathbf{1}_{B(0,r)}}{\vol(B(0,r))} ,
\]
where $\mathbf{1}_A$ denotes the characteristic function of a set
$A$ and $*$ is convolution.  Thus, in particular, the trace of
$\Upsilon _W(z,r)$ is precisely the average area of $W$ in the ball
of radius $r$ and center $z$.}

\medskip

A useful concept in the study of interpolation and sampling for
smooth hypersurfaces is the density of these hypersurfaces.  Let
\[
\varphi _r := \frac {\mathbf{1}_{B(0,r)}* \varphi}{\vol(B(0,r))}  .
\]

\begin{defn}The density of $W$ in the ball of radius $r$ and center $z$
is
\[
D(W,z,r) := \sup \left \{  \frac{\Upsilon _W (z,r)(v,v)}{\ii \di
\dbar \varphi _r (v,v)} \ ;\ v \in T_{\C ^n ,z} - \{0\} \right \}.
\]
The upper density of $W$ is
\[
D^+(W) := \limsup _{r \to \infty} \sup _{z \in \C ^n} D(W,z,r)
\]
and the lower density of $W$ is
\[
D^-(W) := \liminf _{r \to \infty} \inf _{z \in \C ^n} D(W,z,r).
\]
\end{defn}

\rmk {\it Observe that $D(W,z,r) \le (1-c)$ for some $c >0$ if and
only if
\[
\Upsilon _W (z,r) \le (1-c) \ii \di \dbar \varphi_r (z).
\]
On the other hand, a lower bound for $D(W,z,r)$ tells us only that
the largest eigenvalue of the form $\Upsilon _W (z,r) - \ii \di
\dbar \varphi _r (z)$ is uniformly positive.}

\medskip

Our main results can be stated as follows.

\begin{mthm}\label{interp-thm}
Let $W$ be a uniformly flat hypersurface.  If $D^+ (W) < 1$ then
$W$ is an interpolation hypersurface.
\end{mthm}

\begin{mthm}\label{samp-thm}
Let $W$ be a uniformly flat hypersurface.  If $D^-(W)
> 1$ then $W$ is a sampling hypersurface.
\end{mthm}

The hypotheses in Theorems \ref{interp-thm} and \ref{samp-thm}
have a geometric interpretation.  For simplicity, consider the
classical Fock space, which correponds to  $\varphi=|z|^2$.  Then
$\Upsilon _W(z,r)(v,v)$ is the average number of intersections of
the manifold $W$ with a complex line of direction $v$ in the ball
of center $z$ and radius $r$. Thus $D^+ (W) < 1$ means that in any
point $z$ and in any direction $v$ the average number of
intersecting points between the manifold and a complex line in the
direction of $v$ is smaller than some critical value. On the other
hand  $D^-(W) > 1 $ means that for any point $z$ there is a
direction $v$ (which may depend on the point $p$) such that in the
ball of radius $r$ and center $z$ the average number of
intersections between $W$ and the complex line with direction $v$
is bigger than some critical value.

Intuitively speaking, the interpretation of our theorems is  that if we want
$W$ to be interpolating it must be sparse in all points and all directions, but
if we want it to be sampling it must be dense in all points, but only in one
direction for any given point.

\medskip

The interpolation and sampling problems in the generalized Fock
space have been studied previously. In dimension one there is a
full description given in \cite{quimbo} that corresponds to our
Theorem~\ref{interp-thm} and Theorem~\ref{samp-thm}.  In dimension
1 the conditions we require are also necessary.  This was proved
in \cite{quimseep}.  It seems plausible that this is also the case
in higher dimensions, but the question of necessity remains open.

In several complex variables, there have been many partial and related
results. See for instance \cite{BerTay82}, \cite{Berndtsson83} or
\cite{Demailly82}. In these works hypotheses are placed on the  function
$T\in \mathcal O(\C^n)$ defining $W=Z(T)$ in order that $W$ be interpolating
in the sense of Berenstein-Taylor, that is to say, any holomorphic function
$h$ defined on $W$ and satisfying a growth condition
\[
|h(z)|\le C \exp (C\varphi(z))
\]
can be extended to an entire function satisfying similar bounds (perhaps with a
different constant). For instance a result can be found in \cite{BerTay82}
stating that $W$ is interpolating in this sense if
\[
\| \partial T\| \ge C \exp(-C\varphi)\qquad \text{on $W$}.
\]
Our results do not involve the defining function $T$, appealing instead
directly to the current of intergration defining $W$. In this sense our
results are more geometric in nature.

The organization of the paper is as follows.  In Section
\ref{unif-flat-section} we define and discuss the notion of
uniform flatness.  In Section \ref{sing-section} we define a
non-positive function that is singular along the variety $W$.  As
in \cite{quimbo}, this function is used to modify the weight of
the Bargmann-Fock space in order to apply the
H\"ormander-Bombieri-Skoda technique in the proof of
Theorem~\ref{interp-thm}.  A central point is the use of the
Newton potential in the construction.  In Section
\ref{interp-section} we prove Theorem~\ref{interp-thm}. We begin
with the $L^2$ case.  Our approach is to first extend the
candidate function to small neighborhoods, and then to patch
together these local extensions using the solution of a Cousin I
problem with $L^p$ bounds.  To pass to $L^p$ we use results of
Berndtsson on $L^p$ bounds for minimal $L^2$ solutions of $\dbar$.
A second proof is possible when $p=2$, using the method of the
Ohsawa-Takegoshi extension theorem.  This proof is only mentioned
and briefly sketched here. For the details of this approach in the
case of the Bergman ball, the reader is referred to \cite{fv}.  In
Section~\ref{samp-section} we prove Theorem~\ref{samp-thm}.
Finally, in Section~\ref{appl-section} we give a simple
application of our results to improve on known sufficient
conditions for sequences to be interpolating or sampling in $\C
^n,\ n \ge 2$.

\medskip
\noi \textbf{Acknowledgement.}  The authors would like to thank
Tamas Forgacs, Jeff McNeal and Yum-Tong Siu for stimulating and
useful discussions. Some of this work was done while the first
author was visiting the University of Wisconsin, the second author
was visiting the University of Michigan, and the third author was
visiting Harvard University and the University of Michigan.  The
authors wish to thank these institutions for their generous
hospitality.

\section{Uniform flatness} \label{unif-flat-section}

We shall be interested in smooth hypersurfaces $W$ satisfying the
following assumption.
\begin{enumerate}
\item[(F1)]  There is a positive constant $\varepsilon _o$ such that the
Euclidean neighborhood
\[
N_{\varepsilon _o} (W) := \{ z \in \C ^n\ ;\ \operatorname{dist} (z,W) < \varepsilon_o \}
\]
is a tubular neighborhood of $W$ in $\C ^n$.
\end{enumerate}

\begin{defn}
A smooth hypersurface $W$ satisfying \textrm{(F1)} is said to be
uniformly flat.
\end{defn}

If we want to extend functions to $\C^n$, uniform flatness seems a
reasonable condition;  we don't want points that are very far
apart in $W$ to be very close to each other in the ambient space.
When $n=1$, $W$ is a discrete set, which is uniformly flat if and
only if it is uniformly separated.

For each $z \in W$, denote by $T_{W,z}$ the tangent space
to $W$ at $z$ and by $n_z$ a unit normal to $T_{W,z}$ in the
Euclidean metric $\omega$.  Note that $n_z$ is determined uniquely
up to a unimodular constant.  Write
\[
D_W(z,\varepsilon):= \left (T_{W,z} \cap B(z,\varepsilon)\right ) \times \{
\zeta n_z\ ;\ |\zeta | < \varepsilon \}
\]
for the product of the $\varepsilon$-ball centered at the origin in
$T_{W,z}$, with the $\varepsilon$-disc centered at the origin of
$T_{\C ^n,z}$ and perpendicular to $T_{W,z}$.  We leave it to the
reader to verify the following proposition.

\begin{prop}\label{flat-properties}
Let $W$ be a uniformly flat hypersurface.  Then the following hold.

\begin{enumerate}

\item[{\textrm(A)}] For each $R > 0$ there is a constant $C_R
>0$ such that for all $z \in \C ^n$,
\[
\area(W \cap B(z,R)) \le C_R.
\]
\item[{\textrm(B)}]There are positive constants $\varepsilon _o$ and $C$
such that for all $z \in W$, $W\cap D_W (z,\varepsilon_o)$ is given
as a graph over $T_{W,z} \cap B(z,\varepsilon_o)$ by a function $t =
f(x)$, where $f : T_{W,z} \cap B(z,\varepsilon_o) \to \C$ satisfies
\[
|f(z+x)| \le C |x|^2.
\]
\end{enumerate}
\end{prop}

\noi It is not hard to see that every smooth affine algebraic
hypersurface is uniformly flat.  There are also many non-algebraic
examples.

\section{Singularization of the weight}\label{sing-section}

As is now standard in $L^p$ interpolation problems in several
complex variables, one needs to define a strictly plurisubharmonic
weight similar to $\vp$ with singularities along the divisor $W$.
For the sampling problem, one must smooth out this weight near $W$,
while maintaining global bounds away from $W$.

Our scheme for singularizing the weight follows the method of
\cite{quimbo}:  we add to our weight $\varphi$ a function $s_r$,
called the singularity, to be defined below.

To obtain good properties of the singularity, one needs to use
potential theoretic aspects of the ambient space $\C ^n$.  For our
purposes, the Newton potential plays a key role.  Recall that the
Newton potential is the function
\[
G(z,\zeta) = -c(n)|z-\zeta|^{2-2n},
\]
where
\[
c(n) =\frac{1}{\pi ^n 2^{n}(n-1)},
\]
For each $\zeta \in \C ^n$, this function is harmonic in $\C ^n -\{
\zeta \}$ and has the property that
\[
\int _{\C ^n} \ii \di \dbar  G( \cdot ,\zeta ) \wedge \omega
^{n-1} = 1.
\]
The key feature making our approach possible is that this last
identity involves only the trace of $\ii \di \dbar G$.  It is this
fact precisely that links the fundamental solution of $\Delta$ to
holomorphic functions on hypersurfaces.

\subsection*{The singularity}

Consider the function
\[
\Gamma _r (z,\zeta) := \left ( G(z,\zeta)  - \frac{1}{\vol(B(z,r))}
\int _{B(z,r)} G(\zeta ,x) \omega ^n (x) \right ).
\]
Since $G(z,\zeta)$ is harmonic in each variable separately when
$|z-\zeta|>0$, one sees immediately that $\Gamma$ is supported on the
neighborhood $|z-\zeta|\le r$ of the diagonal in $\C ^n \times \C
^n$.  We define the singularity
\begin{eqnarray*}
s_r(z) &:=& \int _{\C ^n} \Gamma _r (z,\zeta) \omega ^{n-1} (\zeta)
\wedge \left ( \ii \di \dbar \log |T| \right ) (\zeta)\\
&=&  \int _{B(z,r)} \Gamma _r (z,\zeta) \omega ^{n-1} (\zeta)
\wedge \left ( \ii \di \dbar \log |T| \right ) (\zeta).
\end{eqnarray*}
By the Lelong-Poincar\'e identity, we have
\begin{eqnarray}\label{post-pl}
s_r(z) = \pi \int _{W _{z,r}}\iifix G(z,\zeta) \omega ^{n-1}
(\zeta) -\frac{\pi}{V(r)} \int _{W_{z,r}}\ifix \left ( \int
_{B(z,r)} \iifix G(\zeta ,x) \omega ^n (x) \right ) \omega ^{n-1}
(\zeta),
\end{eqnarray}
where
\[
W_{z,r} = W \cap B(z,r) \qquad \text{and} \qquad V(r) = \int _{B(z,r)}\omega
^n.
\]
We have the following lemma.
\begin{lem}\label{sing-interp}
The function $s_r$ has the following properties.
\begin{enumerate}
\item\label{positivity} It is non-positive.

\item\label{boundedness} For each $r,\varepsilon > 0$ there is a constant
$C_{r,\varepsilon}$ such
that if $\dist(z, W) \ge \varepsilon$, then $s_r(z) \ge -
C_{r,\varepsilon}$.

\item\label{singular} The function $e^{-2s_r}$  is not locally integrable
at any point of $W$.
\end{enumerate}
\end{lem}

\begin{proof}
By the sub-mean value property for subharmonic functions, $\Gamma _r
\le 0$, from which (\ref{positivity}) follows.  Next, we verify that
there is a constant $D_r$ such that for all $\zeta \in B(z,r)$,
\[
- \frac{1}{\vol(B(z,r))} \int _{B(z,r)}
G(\zeta ,x) \omega ^n (x) \le D_r.
\]
For this, it suffices to bound the integral
\[
I_r(z,\zeta) := \int _{B(z,r+1)}- G(\zeta , x) \omega
^n(x)
\]
Letting $\rho = r - |\zeta|$, we have
\[
I_r(z,\zeta) =  \int _{B(\zeta , \rho+1)} -G(\zeta , x) \omega
^n(x) + \int _{B(z,r)-B(\zeta,\rho +1) }-G(\zeta , x) \omega
^n(x).
\]
Now,
\begin{eqnarray*}
 \int _{B(\zeta , \rho+1)} -G(\zeta , x) \omega ^n(x)
&=& c(n) \int _{B(0,\rho+1)}|y|^{-2n +2} \omega ^n(y)\\
&=& 2 \tilde c(n) \int _0 ^{\rho +1} t dt  \\
&=& \tilde c(n) (\rho +1) ^2 \le \tilde c(n) (r+1)^2.
\end{eqnarray*}
On the other hand,
\begin{eqnarray*}
\int _{B(z,r)-B(\zeta,\rho +1) }\iifix \iifix \iifix \iifix
-G(\zeta , x) \omega ^n(x) &\le & c(n)\int _{B(z,r+1)}\iifix \iifix
\omega ^n,
\end{eqnarray*}
which demonstrates the bound for $I_r$.

If we now look at $z$ such that $|z-\zeta| \ge \varepsilon$ for all $\zeta \in
W$, (\ref{boundedness}) follows from the above bound for $I_r$ together with
formula \eqref{post-pl}.

To prove (\ref{singular}), we study the singularity of the function
\begin{eqnarray}\label{sing-part}
\pi \int _{W \cap B(z,r)} G(z,\zeta) \omega ^{n-1}(\zeta)
\end{eqnarray}
in the neighborhood of a point $w_o \in W$.  Fix Euclidean coordinates
$t, x^1,\ldots,x^{n-1}$ at $w_o$, such that $w_o$ is the origin of these
coordinates, $dt = 0$ defines $T_{W,w_o}$ and
\[
T_{W,w_o} = \operatorname{span}
\left \{ \left . \frac{\di }{\di x^1}\right |_{w_o} ,
\ldots ,  \left . \frac{\di }{\di x^{n-1}}\right |_{w_o} \right \}.
\]
It follows that there are local coordinates $\zeta ^1,\ldots,\zeta
^{n-1}$ on $W$ near $w_o$ that are of the form
\[
\zeta ^i = x^i + O(|x|^2 +|t|^2), \qquad 1 \le i \le n-1.
\]
Moreover, the hypersurface $W$ is cut out by a holomorphic
function of the form
\[
t + O(|t|^2 +|x|^2).
\]
Thus the singularity of the integral \eqref{sing-part} is the same
as that of
\[
\int _{|x|\le c} \frac{-\pi c(n)}{(|t|^2 + |x|^2)^{n-1}} \omega
^{n-1} (x) = \log |t| + O(1).
\]
This completes the proof.
\end{proof}

The proof of (\ref{singular}) in Lemma~\ref{sing-interp} also follows from the
uniform flatness of $W$ and the following formula for $s_r$.

\begin{prop}\label{t-rep}
Let $T \in \co (\C ^n)$ be a holomorphic function such that $W = \{
T=0\}$ and $dT$ is nowhere zero on $W$.  Then
\[
s_r(z) = \log |T(z)| - \frac{1}{V(r)}
\int _{B(z,r)} \log |T(\zeta)| \omega ^n
(\zeta),
\]
and thus
\[
i\ddbar s_r = \Theta_W -\Theta_W * \frac{\mathbf{1}_{B(0,r)}}{\vol(B(0,r))}.
\]
\end{prop}

\begin{proof}
Let $\alpha : [0,\infty ) \to [0,1]$ be a smooth compactly
supported function which is identically 1 on $[0,1]$.  Then for $R >>
r$, we have
\begin{eqnarray*}
s_r(z) &=& \int _{\C ^n} \ii \di \dbar \log |T(\zeta)| \wedge \left (
\Gamma (z,\zeta) \omega ^{n-1}(\zeta )\right ) \\
&=&  \int _{\C ^n} \alpha (R^{-2} |\zeta -z |^2)  \ii \di \dbar \log
|T(\zeta)| \wedge \left ( \Gamma (z,\zeta) \omega ^{n-1}(\zeta )\right ) ,
\end{eqnarray*}
where the second equality follows from the fact that $\Gamma(\cdot , z)$ is
supported on $B(z,r)$.  Integrating by parts and letting $R \to
\infty$, we have
\begin{eqnarray*}
s_r(z) &=& \ifix \left . \int _{B(z,r)} \log |T(\zeta)| \wedge \right  \{
(\ii \di \dbar )_{\zeta} G(z,\zeta) \wedge \omega ^{n-1}(\zeta ) \\
&&  \left . - \left ( \frac{1}{V(r)} \int _{B(z,r)}
\left [ (\ii \di \dbar )_{\zeta} G(x,\zeta) \wedge \omega ^{n-1}
(\zeta) \right ] \omega ^n(x) \right )\right \} \\
&=& \ifix \log |T(z)| -
\left \{ \int _{B(z,r)} \iifix \log |T(\zeta )|
\left ( \frac{1}{4V(r)} \int _{B(z,r)} \iifix
\Delta _{\zeta} G(x,\zeta) \omega ^n(x) \right )
\omega ^{n} (\zeta) \right \} \\
&=&\ifix \log |T(z)| -
\left \{ \int _{B(z,r)} \iifix \log |T(\zeta )|
\left ( \frac{1}{4V(r)} \int _{B(z,r)} \iifix
\Delta _{x} G(x,\zeta) \omega ^n(x) \right )
\omega ^{n} (\zeta) \right \} \\
&=& \ifix \log |T(z)| - \frac{1}{V(r)} \int _{B(z,r)} \log
|T| \omega ^n,
\end{eqnarray*}
as desired.
\end{proof}

\section{Interpolation:  The proof of Theorem~\ref{interp-thm}}
\label{interp-section}
Since we assume that $i\ddbar \varphi \simeq \omega$, it follows that
$|\varphi_r-\varphi|\le C$ and therefore the spaces $\BF_\varphi^p$ and
$\BF_{\varphi_r}^p$ are the same space with equivalent norms. The same happens
with $\lf_\varphi^p$ and $\lf_{\varphi_r}^p$. Therefore we may assume without
loss of generality that in the definition of the densities and thus in the
hypothesis of the theorems we have replaced $\varphi_r$ by $\varphi$.

\subsection{The Cousin I approach}

\subsubsection*{Local extension}

Let $\varphi$ be a plurisubharmonic function in $\C^n$ with $\ii \di
\dbar \varphi \le M \omega$ for some $M >0$. Let $\Omega$ be a bounded
domain in $\C^n$ containing the origin and denote by $H$ the
hyperplane $z_n=0$.  Define
\[
H_{\Omega} := P_H (\Omega),
\]
where $P_H$ denotes orthogonal projection onto $H$.

\begin{prop} Assume that $P_H(\Omega)=\Omega \cap \{z_n=0\}$.
There is a constant $C>0$, depending only on $M$ and on the
diameter of $\Omega$, such that for any holomorphic function $f\in
\co (H_{\Omega})$ there is a function $F\in \co (\Omega)$ such
that $F|H\cap \Omega = f$ and
\[
\int_\Omega |F|^{p}e^{-p\varphi} \omega ^n \le C^p \int_{H_{\Omega}}
|f|^{p}e^{-p\varphi} \omega ^{n-1},\quad 0 < p \le \infty,
\]
provided that the right hand side is finite.
\end{prop}
\noi (When $p=\infty$, the integrals should be replaced by suprema
and $C^p$ by $C$.)

\begin{proof}
Let $D_{\Omega}$ be a disk such that
\[
\Omega \subset H_{\Omega} \times D_{\Omega} \subset B(0,R),
\]
where $R = \operatorname{diam} (\Omega)$.  In $B(0,R)$ there is a function
$u$ such that
\begin{enumerate}
\item $u$ is bounded in $B(0,R)$ by a constant depending only on
$M$ and $R$, and

\item
\[
\ii \di \dbar u = \ii \di \dbar \varphi.
\]
\end{enumerate}
\noi (For the proof, see \cite[Lemma 6]{l-01}.)  Define $h= \varphi-
u$. Since $h$ is pluriharmonic there is a function $H \in \co
(B(0,R))$ such that $\Re\ H =h$.  Writing $z = (z',z^n) \in \C
^{n-1} \times \C$, we let
\[
F(z',z^n) := f(z') e^{H(z',z^n) - H(z',0)}.
\]
Then $F \in \co (H_{\Omega} \times D_{\Omega})$, $F|H_{\Omega}=
f$, and we have
\begin{eqnarray*}
|F(z)|^p e^{-p\varphi(z)}&=& |f(z')|^p  \exp(\varphi(z',z_n)-\varphi(z',0))
e^{-p\varphi(z)}\\
&=& |f(z')|^p \exp(-pu(z)+pu(z',0)-p\varphi(z',0))\\
&\le& K |f(z')|^p e^{-p\varphi(z',0)}.
\end{eqnarray*}
The result follows.
\end{proof}

\medskip

\noi As a corollary, we have the following lemma.

\begin{lem}\label{grad-t-extension}
Let $W$ be uniformly flat, let $0<p \le \infty$ and let $\varepsilon <
\varepsilon
_o/2$, where $\varepsilon _o$ is as in \ref{flat-properties}-\textrm{(B)}.
Then there is a
constant $C_p>0$ depending only on $\varepsilon _o$ $p$ and $M$, such that the
following holds. For each $w \in W$ and $f \in \co (B(w,2 \varepsilon)\cap
W)$ there is a function $F \in \co (B(w,\varepsilon))$ such that
\[
F|B(w,\varepsilon)\cap W = f \quad \text{and} \quad \int _{B(w,\varepsilon)}
\!\!\!\! \!\!\!\! |F|^{p} e^{-p\varphi} \omega ^n \le C \int
_{W\cap B(w,2\varepsilon)}\!\!\!\! \!\!\!\! \!\!\!\!\!\! |f|^{p}
e^{-p\varphi}\omega ^{n-1}.
\]
If $p=\infty$, then the integrals should replaced by suprema.
\end{lem}

\subsubsection*{Local holomorphic functions with good estimates}

\begin{lem}\label{quimbo-trick}
Let $\varphi$ be a function in the unit disk $\D$ such that
\[
c \le \Delta \varphi \le \frac{1}{c}.
\]
Then there exist a constant $C>0$ and a holomorphic function $H\in
\co (\D)$ such that $H(0) = 0$ and
\[
| \Re H - \varphi + \varphi (0)| \le C.
\]
Moreover, if $\varphi$ depends on a parameter in such a way that the
bound on $\Delta \varphi$ is independent of the parameter, then $H$
can be taken to depend on this parameter in such a way that $C$
does not.
\end{lem}
\noi The proof of this lemma, by now well known, can be found in
\cite{quimbo}.

\subsubsection*{Construction of the interpolating function}

We fix $f \in \lf^p_{\varphi}(W)$ and $\varepsilon < \varepsilon _0 /2$, where
$\varepsilon _0$ is as in \ref{flat-properties}-(B).
Take a sequence of distinct points $\{
w_j \ ;\ j=1,2,\ldots\} \subset W$ such that
\[
N_{\varepsilon}(W) \subset \bigcup_{i=1}^\infty
\left \{B\left (w_i,\tfrac{3}{2}\varepsilon\right )\right \}_{w_i\in W}
\]
and each point of $N_{\varepsilon}(W)$ is contained in at most a fixed,
finite number of the sets
\[
B(w_j,2\varepsilon).
\]
(We say that the cover is uniformly locally finite.)  For
convenience of notation we write $B_i = B\left (w_i,\tfrac{3}{2}
\varepsilon\right )$. We add to the cover $\{ B_i \}_{i \ge 1}$
another open set $B_0=\C^n\setminus N_{\frac{1}{2}\varepsilon}(W)$.
Thus $\{B_j\ ;\ j \ge 0 \}$ is a uniformly locally finite open cover
of $\C^n$. Let $\{ \phi _i\}_{i \ge 0}$ be a partition of unity
subordinate to the cover $\{ B_i\}$, i.e., $0\le\phi_i\le 1$,
$\operatorname{supp} \phi_i \in B_i$ and $\sum_i \phi_i \equiv 1$.
Moreover we can assume that $\sum _i \|d\phi_i\| \le C$.

Let $F_i$ denote the extension to $B_i$ of $f|W\cap B(w_i,2\varepsilon)$
given by Lemma~\ref{grad-t-extension}, and set $F_0 \equiv 0$.
Since the covering $\{ B_i\}$ is uniformly locally finite, we have
\[
\int_{\C^n} \sum_{i} \chi _i |F_i|^pe^{-p\varphi}\omega ^n \lesssim
\int _W |f|^p e^{-p\varphi} \omega ^{n-1},
\]
where $\chi _i$ denotes the characteristic function of $B_i$ and, as
usual, the symbol $\lesssim$ means that the left hand side is
bounded above by a {\it universal} constant times the right hand
side. We want to patch together the extensions $F_i$ and construct a
single holomorphic extension $F$ of $f$ whose norm remains under
control. In the standard language of several complex variables, we
want to solve a Cousin I problem with $L^p$ bounds. The setup of the
problem is as follows. For any pair of indices $i,j \ge 0$ we define
a function $G_{ij}$ in $B_{ij} := B_i\cap B_j$ by
\[
G_{ij}=F_i-F_j.
\]
Observe that
\[
G_{ij}|W \cap B_{ij} \equiv 0 \quad \text{and} \quad
G_{ij}+G_{jk}+G_{ki}\equiv 0 \text{ in }  B_i\cap B_j \cap B_k.
\]
Finally
\[
\int_{\C^n} \sum_{i,j} \chi _i \chi _j |G_{ij}|^pe^{-p\varphi} \omega
^n\lesssim \int _W |f|^p e^{-p\varphi} \omega ^{n-1}.
\]
We seek $G_i\in \co (B_i)$ such that $G_{ij}=G_i-G_j  \text{ in }
B_{ij}$, $G_i|W\cap B_i \equiv 0$ and
\[
\int_{\C^n} \sum_{i} \chi _i |G_i|^pe^{-p\varphi} \omega ^n  \lesssim
\int _W |f|^p e^{-p\varphi} \omega ^{n-1}.
\]
If we find such functions $G_i$, then the function $F$ defined
by
\[
F(x) =F_i(x)-G_i(x)\quad x\in B_i
\]
is an entire function.  (It is well defined because $F_i-F_j = G_i
- G_j$ on $B_{ij}$.)  Moreover we have
\[
F |_W=f \quad \text{ and } \quad \int_{\C^n} |F|^p
e^{-p\varphi}\omega ^n \lesssim \int _W |f|^p e^{-p\varphi}\omega ^{n-1}.
\]

We define $\tilde G_i \in \mathcal{C}^{\infty} (B_i)$ by $\tilde
G_i=\sum_j \phi_j G_{ij}$. These functions have all the properties
we seek, except they are not holomorphic.  We shall now correct
the functions $\tilde G_i$ by adding to each of them a single,
globally defined function.

To this end, note that in $B_{ij}$ we have $\dbar \tilde G_i =
\dbar \tilde G_j$.  Thus there is a well defined $\dbar$-closed
$(0,1)$-form $h$ such that
\[
h=\dbar \tilde G_i \text{ in }\ B_i.
\]
Moreover, observe that
\[
\|h\| \le \sum_{ij}\left \| \dbar \phi_i \right \|
\cdot |G_{ij}|.
\]
\begin{lem}\label{h-L2-est}
One has the estimate
\[
\int _{\C ^n} \|h\|^p e^{-p(\varphi + s_r)} \omega ^n \le C \int _W
|f|^p e^{-p\varphi} \omega ^{n-1} .
\]
\end{lem}

\begin{proof}
Recall that if $\psi$ is a weight function on the unit disk
$\D$ such that $\Delta \psi \le K$, then there is a constant $C$
such that for any $f\in \co (\D)$,
\[
\int_{|z|<1} |f|^p e^{-p\psi} \le K \int_{1/2<|z|<1} |f|^p
e^{-p\psi}.
\]
Indeed, the inequality is elementary in the case $\psi\equiv 0$.
Since the Laplacian of $\psi$ is bounded, there is by
Lemma~\ref{quimbo-trick} a non-vanishing holomorphic function $g$, such
that $|g|\simeq e^{\psi}$.  Thus we obtain
\[
\int_{|z|<1} |f|^p e^{-p\psi} \simeq \int_{|z|<1} |f/g|^p \le K
\int_{1/2<|z|<1} |f/g|^p \simeq \int_{1/2<|z|<1} |f|^p e^{-p\psi}.
\]
With this one variable fact it is possible to prove that
\begin{eqnarray*}
\int_{B_i\cap B_j} |G_{ij}|^p e^{-p(\varphi+s_r)} &\lesssim&
\int_{(B_i\cap B_j)\setminus N_{\frac{1}{2}\varepsilon}(W)} |G_{ij}|^p
e^{-p(\varphi+s_r)}\\
&\simeq&
\int_{(B_i\cap B_j)\setminus N_{\frac{1}{2}\varepsilon}(W)} |G_{ij}|^p
e^{-p\varphi}\\
&\lesssim& \int_{(B_i\cup B_j)\cap W} |f|^p e^{-p\varphi}.
\end{eqnarray*}
Only the first inequality is non-trivial.  To see how it follows,
let $T$ be any entire function that vanishes precisely on $W$ such
that $dT$ does not vanish on $W$.  Then by Proposition~\ref{t-rep},
\[
s_r=\log|T|-\frac{\mathbf{1}_{B(0,r)}}{\vol(B(0,r))} * \log|T|.
\]
Therefore
\[
|G_{ij}|^p e^{-p(\varphi + s_r)}\simeq |G_{ij}/T|^p e^{-p\psi_r},
\]
where
\[
\psi_r =\varphi - \frac{\mathbf{1}_{B(0,r)}}{\vol(B(0,r))} * \log |T|.
\]
It follows by the density hypothesis that
\[
\ii \di \dbar \psi_r\simeq Id.
\]
Since the function $G_{ij}/T$ is holomorphic in $B_i\cap B_j$, we
may apply the one dimensional result above.  Let $U=B_{ij}\cap W$.
Then $B_{ij} \simeq U\times D(0,\varepsilon)$. We integrate along
the slices and apply the one-dimensional result in each disk.
\end{proof}

\medskip

By the density hypothesis, one has the inequality $ \ii \di \dbar
(\varphi+s_r) \ge c \omega
> 0.$
We will deal first with the case $p=2$.  It follows from H\"ormander's Theorem
that there is a function $u$ such that
\[
\dbar u = h \quad \text{ and } \quad \int _{\C ^n} |u|^2 e^{-2(\varphi
+s_r)} \omega ^n \le C \int _{W} |f|^2 e^{-2\varphi}\omega ^{n-1}.
\]
Moreover, the local non-integrability of $e^{-2(\varphi +s_r)}$ on $W$
guarantees that $u|W\equiv 0$. Finally, since $\varphi \ge \varphi +s_r$,
we have that
\[
\int _{\C ^n} |u|^2 e^{-2\varphi} \omega ^n \le \int _{\C ^n} |u|^2
e^{-2(\varphi +s_r)} \omega ^n.
\]
It follows that the holomorphic functions $G_i= \tilde G_i-u$ have
the desired properties.

Next we treat the case $p \in [1,2)$.   Let us denote $\xi
=\varphi+s_r$.  Since $h$ is supported away from the singularity
of $\xi$, a look at the definition of $h$ (in particular, it is
constructed from certain holomorphic data and cutoff functions)
shows that, since $h \in L^p(e^{-\xi})$, $h \in L^{\infty}
(e^{-\xi})$.  It follows that $h \in L^2 (e^{-\xi})$.  Let $u$ be
the function of minimal norm in $L^2 (e^{-\xi})$ satisfying $\dbar
u = h$. Then a theorem of Berndtsson
\cite{Berndtsson97,Berndtsson01} states that $u$ satisfies
\[
\| u e^{-\varphi}\|_{L^p}\le C_p \|he^{-\xi}\|_{L^p},\quad p\in[1,\infty],
\]
provided the right hand side is finite (which in the case at hand applies
for $p \in [1,2]$), and  $\ii\ddbar \varphi \simeq \omega$, $\ii\ddbar\xi \ge c
\omega$ and $s_r\le 0$, as is indeed the case here.  (We point out that the
constants $C_p$ in Berndtsson's Theorem depend only on $p$ and on
the upper and lower bounds for $\ii \di \dbar \vp$.)
This gives the right bounds for the solution.
Moreover, since
$\|he^{-\xi}\|_{L^2} < +\infty$, H\"ormander's Theorem and
the minimality of $u$ tell us that $\|ue^{-\xi}\|_{L ^2}<+\infty$.
Thus again $u|W\equiv 0$.

Finally, we come to the case $p\in (2,\infty]$.  Here we must be a
little more careful. Assume first that $h \in L^2(e^{-\xi}) \cap
L^p(e^{-\xi})$.  Let $u$ be the function of minimal norm in $L^2
(e^{-\xi})$ such that $\dbar u = h$.  Then again by Berndtsson's
Theorem $u$ satisfies
\[
\| u e^{-\varphi}\|_{L^p}\le C_p \|he^{-\xi}\|_{L^p},\quad p\in[2,\infty].
\]
This again gives the desired bounds.  Moreover, if
$\|he^{-\xi}\|_{L^2}$ is finite then by H\"ormander's Theorem and
the minimality of $u$, we have $\|ue^{-\xi}\|_{L ^2}<+\infty$.
Thus again $u|W\equiv 0$.

This proves the result for $h \in L^2(e^{-\xi}) \cap L^p(e^{-\xi})$.
To pass to the general case, instead of approximating $h$ we modify
the weight $\xi$.  To this end, take any sequence $\varepsilon _j
\to 0$.  Since $h$ is identically zero on a neighborhood of $W$ and
$he^{-\xi} \in L^p$, we have $he^{-\vp} \in L^p$.  Thus once again
$he^{-\vp} \in L^{\infty}$, and by the support of $h$ we have
$he^{-\xi} \in L^p$. It follows that for all $j$,
$he^{-\xi-\varepsilon_j\|z\|^2}\in L^2$. As before, the solution
$u_j$ to $\dbar u_j = h$ with minimal norm in
$L^2(e^{-\xi-\varepsilon_j\|z\|^2})$ vanishes on $W$ and, by
Berndtsson's Theorem, satisfies
\[
\| u_j e^{-\vp-\varepsilon_j\|z\|^2}\|_{L^p}\le C_p
\|he^{-\xi-\varepsilon_j\|z\|^2}\|_{L^p},\quad p\in[1,\infty],
\]
where the constants $C_p$ are independent of $j$. It follows that $u
_j \to u \in L^p (e^{-\vp})$. Thus we can construct holomorphic
functions $F^j$ that extend $f$ and satisfy the estimates
\[
\int |F^j|^p e^{-p\vp-p\varepsilon_j\|z\|^2} \le C_p
\int_{W} |f|^p e^{-p\vp}.
\]
By a normal family argument we can take a subsequence $F^j$
converging to $F\in L^p(e^{-\vp})$.
The convergence is unifom over
compacts and thus $F$ extends $f$.

\subsection{Remarks on the twisted $\dbar$ approach}

In this section we outline the ideas behind a proof of
Theorem~\ref{interp-thm} in the case $p=2$
using the method of the twisted $\dbar$ equation.

The idea behind the twisted $\dbar$ approach is to replace the
usual $\dbar$ complex
\[
A ^2 _{0,0}(\Omega , \psi) {\buildrel \dbar \over \longrightarrow}
A ^2 _{0,1}(\Omega , \psi) {\buildrel \dbar \over \longrightarrow}
A ^2 _{0,2}(\Omega , \psi)
\]
by a complex
\[
A ^2 _{0,0}(\Omega , \psi) {\buildrel T \over \longrightarrow} A
^2 _{0,1}(\Omega , \psi) {\buildrel S \over \longrightarrow} A ^2
_{0,2}(\Omega , \psi),
\]
where the two operators $T$ and $S$ are defined by
\[
Tu = \dbar ((\sqrt{\tau +A})u) \quad \text{and} \quad Su = \sqrt
\tau (\dbar u).
\]

If the domain $\Omega \subset \subset \C ^n$ is smoothly bounded and
pseudoconvex, then clever manipulation of the usual
Bochner-Kodaira identity can be used to show that for any
$(0,1)$-form $u$ in the domains of $T^*$ and $S$, a \emph{twisted}
Bochner-Kodaira inequality holds:
\[
\|T^*u\|_{\psi} ^2 + \|Su\|_{\psi} ^2 \ge \int _{\Omega}  \left (
\tau \ii \di \dbar \psi (u,u) - \ii \di \dbar \tau (u,u) - \frac{
\left |\di\tau (u)\right |^2}{A} \right ) e^{-\psi} \omega ^n.
\]
(Actually, the best way to obtain this identity is by the method of McNeal-Siu
\cite{m-96,siu-96} of twisting the weights in the usual Bochner-Kodaira
identity.)
By choosing
\[
\psi = \kappa + s_r, \quad \tau = a + \log a \quad \text{and} \quad
A = (1+a)^2,
\]
where
\[
a = 1 +\log (1+\varepsilon ^2) - \log ( e^{s_r} +\varepsilon ^2),
\]
one can deduce from the twisted Bochner-Kodaira inequality an a
priori identity which can be used to solve the equation $Th =
\alpha$ with estimates
\[
\int _{\Omega} |h|^2 e^{-\psi} \omega ^n \le C
\]
whenever $\alpha$ is an $S$-closed $(0,1)$-form such that for all
$u$ with compact support in $\Omega$
\[
|(u,\alpha)|^2 \le C (\|T^*u\|^2 +\|Su\|^2).
\]
The choice of
\[
\alpha = \dbar \chi \tilde f,
\]
where $\tilde f$ is any holomorphic extension of $f \in \co (W)$
to $\C ^n$ and $\chi$ is an appropriate cut-off function,
produces a holomorphic function
\[
F = \chi \tilde f - \sqrt{\tau +A} h.
\]
One estimates this function and passes to the limit as $\Omega \to
\C ^n$, using the Cauchy estimates to pass from $L^2$ convergence
to locally uniform convergence.

As already mentioned, the details of this approach will not be
fully carried out here.  For an adaptation in the case of the
Bergman ball, see \cite{fv}.

\section{Sampling}\label{samp-section}

In this section we prove Theorem~\ref{samp-thm}. As in
section~\ref{interp-thm}, we replace $\varphi$ by
$$
\varphi _r:= \frac{\mathbf{1}_{B(0,r)}*\vp}{\vol(B(0,r))}
$$
in the definition of the density and thus in the hypothesis of
Theorem~\ref{samp-thm}.

\subsection*{Restrictions and the upper sampling inequality}
\begin{prop}\label{restriction}
If $W$ is a uniformly flat hypersurface, then there is a constant
$C>0$ such that for all $F \in \BF_{\varphi}^p (N_{\varepsilon}(W))$
one has
\[
C \varepsilon ^2 \int _W |F|^p e^{-p\varphi} \omega ^{n-1} \le \int
_{N_{\varepsilon}(W)} |F|^p e^{-p\varphi} \omega ^n.
\]
\end{prop}

\begin{proof}
By our hypotheses, $N_{\varepsilon}(W)$ is foliated by analytic
disks, each of which is transverse to $W$ as well as to the
boundary of $N_{\varepsilon} (W)$, and  meets $W$ at a single
point.  For a given $x \in W$, we denote by $\cl _x$ the disk
passing through $x$, and by $\lambda _x : \D \to \cl _x$ the
(unique up to precomposition by a rotation) holomorphic
parameterization of $\cl _x$ by the unit disk, sending $0$ to $x$.

We begin with the following claim:
\[
\area( \cl _x) = \int _{\D} \lambda _ x^* \omega = \int
_{\cl _x} \omega \ge 2\pi \varepsilon ^2.
\]
To see this, let $p \in \cl _x$ be a point on the boundary of
$N_{\varepsilon}(W)$ that is of minimal distance to $x$.  By definition of
$N_{\varepsilon}(W)$,  the distance from $p$ to $x$ is at least $\varepsilon$.
Let $\ell$ be the complex affine line in $\C ^n$ containing $x$
and $p$. Then it follows from our choice of $p$ and from the
maximum principle that the projection of $\cl _x$ onto $\ell$
contains the Euclidean disk in $\ell$ of center $x$ and radius
$\varepsilon$.  Thus the area of $\cl _x$ is at least $2\pi \varepsilon ^2$.

Making use of the diffeomorphism
\[
W \times \D \ni (x,t) \mapsto \lambda _x(t) \in N_{\varepsilon}(W),
\]
which is holomorphic in the second variable, we work on the
product $W \times \D$.

Let $H(x,t)$ be the function, holomorphic in $t$, given by
Lemma~\ref{quimbo-trick}. That is to say,
\[
H(x,0)=0 \qquad \text{and} \qquad \left | \Re(H(x,t)) -\varphi (x,0)
+ \varphi (x,t) \right | \le C
\]
for some positive constant $C$, since we have assumed that $\ii
\di \dbar \varphi$ is bounded above.  We then have
\begin{eqnarray*}
\varepsilon ^2|F (x,0) |^pe^{-p\varphi (x,0)}
&=& \varepsilon ^2 \left |F (x,0) e^{H(x,0)}\right |^p e^{-p\varphi (x,0)}\\
&\le & \frac{1}{2\pi} \int _{\D} \left |F (x,t)e^{H(x,t)}\right
|^p e^{-p\varphi (x,0)} \lambda _x ^* \omega  \\
& \le & C \int _{\D} \left |F (x,t)\right |^p e^{-p\varphi
(x,t)} \lambda _x ^* \omega \\
&=& C \int _{\cl _x} |F |^p e^{-p\varphi} \omega .
\end{eqnarray*}
Integration over $W$ then yields
\begin{eqnarray*}
\varepsilon ^2 \int _W |F|^p e^{-p\varphi} \omega ^{n-1} &\le& C \int
_{N_{\varepsilon}(W)} |F|^pe^{-p\varphi} \wedge \omega ^{n},
\end{eqnarray*}
and the proof is complete.
\end{proof}

\begin{cor}\label{upper-sampling-ineq}
If $W$ satisfies {\rm (F1)} then there is a constant $M > 1$ such
that for every $F \in \BF_{\varphi}^p (\C ^n)$,
\[
\int _W |F|^p e^{-p\varphi} \omega ^{n-1} \le M \int _{\C ^n}
|F|^p e^{-p\varphi} \omega ^n.
\]
\end{cor}

\subsection*{The proof of Theorem \ref{samp-thm}}\ 

\ms

\noi The proof will be an almost immediate application of the following 
sequence of definitions and lemmas.

\begin{defn}
A sequence of complex hypersurfaces $W_n$ is said to converge
weakly to another complex hypersurface $W$ if the corresponding
currents of integration $\Theta_{W_n}$ converge to $\Theta _W$ in
the sense of currents.
\end{defn}

\begin{lem}
If $W$ is a uniformly flat complex hypersurface, then for any
sequence of translations $\tau_n$, the sequence $W_n=\tau_n(W)$
has a subsequence converging weakly to a uniformly flat complex
hypersurface $V$.  Moreover, $V$ has a tubular neighborhood of at 
least the same thickness as that of $W$.
\end{lem}
\begin{proof}
We denote by $|\Theta_{W_n}|$ the trace of the current $\Theta _{W_n}$.  This
is a positive measure that dominates all the coefficients of $\Theta_{W_n}$. 
By the uniform flatness of $W$ it is clear that for any ball $B$, 
$\sup_n |\Theta_{W_n}|(B)<C$ for some constant $C$ depending only on the 
radius of $B$.  A standard compactness argument produces a  
subsequence that converges to a positive closed current $\theta$. 
It remains to show that the limit current $\theta$ is a current of 
integration on a manifold $V$. This is proved in \cite{bishop64}, again 
under the assumptions that for any fixed ball $B$ the mass $|\Theta
_{W_n}|(B)$ is bounded. Moreover, in this situation the support of 
$\Theta_{W_n}$ converges to $V$ and in any ball the tubular neighborhoods 
of the $W_n\cap B$ converge to a tubular neighborhood of $V\cap B$.
\end{proof}

\begin{defn}
A sequence of plurisubharmonic functions $\vp _n$ is said to
converge weakly to a plurisubharmonic function $\vp$ if the
corresponding currents $\ii\ddbar \vp_n$ converge to
$\ii\ddbar\vp$ in the sense of currents.
\end{defn}

\begin{lem}
If $\vp$ satisfies $\ii\ddbar\vp \simeq \omega$, then for any
sequence of translations $\tau_n$, the sequence $\vp_n=\vp\circ
\tau_n$ has a subsequence converging weakly to a plurisubharmonic
$\psi$ and $\ii \ddbar\psi \simeq \omega$, with the constants in
the estimates $\ii \di \dbar \psi \simeq \omega$ controlled by the 
constants in the estimate $\ii \ddbar \vp = \omega$.
\end{lem}

\begin{proof}
This is proved in dimension 1 in \cite{quimseep}. The same proof
applies {\it mutatis mutandi}, so we content ourselves with but a sketch.
Let $\theta_n=\ii\ddbar \vp_n$.  In view of the hypothesis $\ii\ddbar
\vp \lesssim \omega$,  we see that $|\theta_n|(B(z,R)) \le
C_R^{n}$ where $C_R ^n$ is independent of $z$, and there are
functions $\psi_n$ such that $\ii\ddbar \psi_n=\theta_n$,
$\psi_n(0)=0$ and $\ii\ddbar \psi_n$ is uniformly Lipschitz.  By a
normal family argument we can take a subsequence, still denoted
$\psi_n$, such that $\psi_n\to\psi$ uniformly on compacts, and
$\ii\ddbar\psi_n\to \ii\ddbar\psi$ as currents.
\end{proof}

\begin{defn}
Given a pair $(W,\vp)$ where $W$ is a uniformly flat complex
hypersurface and $\vp\in PSH(\C^n)$ with $\ii \ddbar\vp\simeq
\omega$, we denote by $K^*(W,\vp)$ the collection of all pairs
$(V,\psi)$ for which there is a sequence of translations $\tau_n$
such that $\tau_n(W)$ converge weakly to $V$ and $\vp \circ
\tau_n$ converge weakly to $\psi$.
\end{defn}

\begin{lem}\label{estabilitatdensitat}
If the pair $(W,\vp)$ satisfies $D_\vp^{-}(W)= \alpha$ then all
pairs $(V,\psi)\in K^*(W,\vp)$ satisfy $D_\psi^{-}(V)\ge \alpha$
\end{lem}

\begin{proof}
By hypothesis, for any $z \in \C ^n$ and $\ve > 0$ there exists $r>0$ and $v
\in \C ^n$ of unit norm such that
\[
\int_{B(z,r)} \Theta_W(v,v) \ge (1-\varepsilon)\alpha
\int_{B(z,r)} \ii\ddbar \vp(v,v).
\]
We fix an arbitrary $z\in \C^n$. Take a sequence of translations
$\tau_n$ such that $W_n=\tau _n(W)$ and $\vp_n = \tau _n ^* \vp$
converge to $V$ and $\psi$ respectively. By definition of $D_\vp^-(W)=\alpha$, 
for any $\varepsilon>0$, there is an $r>0$ and unit vectors $v_n$ such that
\[
\int_{B(z,r)} \Theta_{W_n}(v_n,v_n) \ge (1-\varepsilon)\alpha
\int_{B(z,r)} \ii\ddbar \vp_n(v_n,v_n).
\]
By compactness there is a subsequence of the $v_n$ converging to
$v$ with $||v||=1$. By Hurwitz's theorem
$$
\liminf_n \int_{B(z,r)} \Theta_{W_n}(v_n,v_n)\le \int_{B(z,r)}
\Theta_{V}(v,v),
$$
and since $\ii\ddbar \vp \simeq \omega$,
$$
\lim_n \int_{B(z,r)} \ii\ddbar\vp_n(v_n,v_n)= \int_{B(z,r)}
\ii\ddbar\psi(v,v).
$$
\end{proof}
\begin{defn}
The pair $(V,\psi)$ is said to be determining if for any $f\in
\BF_\psi^\infty(\C^n)$, $f|_V=0$ implies that $f \equiv 0$.
\end{defn}

\begin{lem}\label{beurlinglem}
The manifold $W$ is sampling for $\BF^\infty_\vp$ if all pairs
$(V,\psi)\in K^*(W,\vp)$ are determining.
\end{lem}

Lemma \ref{beurlinglem} was essentially proved by Beurling in \cite[pp.
341--365]{beurling}, so we omit the proof.  This is a key result
because it allows us to determine that $W$ is sampling simply by checking
the more easily verified condition that $V$ is determining.  


\begin{lem}\label{jensen}
If $D_\psi^-(V)>1$ then the pair $(V,\psi)$ is determining.
\end{lem}
\begin{proof}
Without loss of generality we assume that $0\notin V$. In order to 
arrive at a contradiction, assume there exists an $F\in \BF^\infty
_\vp$ with $F|V\equiv 0$ and $F(0)=1$.  By hypothesis there is a 
direction $v$ such that the density of $V$ in the direction of $v$ 
is greater than 1. We will work on the line $\ell = \C v$.  Write $f=F|\ell$
and $\phi = \vp |\ell $, and let $\Gamma = V \cap \ell$.  Then $\Gamma$ 
is a uniformly separated sequence with density $>1$ with respect to the
weight $\phi$.  Recall that the one-dimensional lower density is
$$
\liminf _{R \to \infty} \inf _{z\in \ell} \frac{\# (\Gamma \cap
D(z,R))}{\int _{D(z,R)}\Delta \phi}.
$$

By hypothesis $f(0)=1$.  Now, if $n(0,s)$ denotes the number of
zeros of $f$ in $D(0,s)$, then
\begin{equation}\label{densitat}
\liminf_{R\to\infty}\frac{n(0,R)}{\int_{D(0,R)}\Delta \phi }>1.
\end{equation}
Applying Jensen's Formula to $f$, we get
\[
\int_1^R \frac{n(0,s)}s\,ds \le \frac{1}{2\pi}\int_0^{2\pi} \log
|f(Re^{\ii\theta})|\, d\theta.
\]
Since $\log|f(Re^{\ii\theta})|\le \phi (Re^{\ii\theta}) +K$, we obtain
\[
\int_1^R \frac{n(0,s)}s\,ds \le\frac{1}{2\pi}\int_0^{2\pi} \phi 
(Re^{\ii \theta}) \, d\theta + K.
\]
Now, by Green's Theorem we have 
\begin{eqnarray*}
\int _0 ^{2\pi} \phi (Re^{\ii \theta})d\theta 
&=& \int _0 ^R \frac{1}{s} \int _0 ^{2\pi} \left ( s \frac{\di }{\di s}  
\phi (se^{\ii \theta})d\theta \right ) ds \\
&=&  4 \int _0 ^R \frac{\int _{D(0,s)} \Delta \phi}{s} ds,  
\end{eqnarray*}
and thus
\[
\int_1^R \frac{n(0,s)}s\,ds\le\frac 2\pi \int_0^R
\frac{\int_{D(0,s)}(\Delta \phi) }{s}\,ds+ K.
\]
Thus since $\int_{D(0,R)}\Delta\phi \simeq R^2$,
\[
\frac{\int_1^R \frac{n(0,s)}s\,ds}
{\int_0^R\frac{\int_{D(0,s)}(\Delta \phi) }{s}\,ds} \le 1+ K/R^2.
\]
which contradicts \eqref{densitat}.

\end{proof}
\begin{lem}\label{discretitzat}
If $W$ is a uniformly flat sampling hypersurface for
$\BF^\infty_\vp$ then there is a uniformly separated sequence
$\Sigma\subset W$ that is sampling for $\BF^\infty_\vp$.
\end{lem}

\rmk {\it The definition of a sampling sequence is given in
Section~\ref{appl-section} below.}

\begin{proof}
Any set $F$ that is sampling for $\BF^\infty_\vp$ contains a
uniformly separated sampling sequence. This is proved in
\cite[Proposition~19]{l-01}.  (For the 1 dimensional case, see
\cite[Proposition~2]{quimseep}.)
\end{proof}
\begin{lem}\label{inclusions}
Let $1\le p \le \infty$. If $\Sigma$ is a uniformly separated
sampling sequence for $\BF^\infty_{\vp+\varepsilon|z|^2}$ then it
is a sampling sequence for $\BF^p_\vp$.
\end{lem}
\begin{proof}
Denote by $\BF^{\infty,0}_{\vp+\varepsilon|z|^2}$ the closed
subspace of $\BF^\infty_{\vp+\varepsilon|z|^2}$ consisting of
functions $f$ such that
$$
\lim_{z\to\infty} |f|e^{-\vp-\varepsilon|z|^2}=0.
$$
The restriction operator
\[
R:\BF^{\infty,0}_{\vp+\varepsilon|z|^2}\to
\ell_{\vp+\varepsilon|z|^2}^{\infty,0}
\]
sending $f$ to $\{f(\sigma)\}_{\sigma\in\Sigma}$ is a bounded
linear operator. Since $\Sigma$ is sampling, $R$ is onto and has
closed range. Thus R defines an isomorphism between
$\BF^{\infty,0}_{\vp+\varepsilon|z|^2}$ and its image. For any
$z\in\C^n$ the weighted point evaluation 
$$
f \mapsto f(z)e^{-\vp(z)-\ve|z|^2}
$$ 
is bounded on $\BF^{\infty,0}_{\vp+\varepsilon|z|^2}$. 
Thus, for every $z$ there is a sequence $k(z,\sigma)$ such that
\begin{equation}\label{representacio}
f(z)e^{-\vp(z)-\varepsilon|z|^2} = \sum_{\sigma\in\Sigma}
k(z,\sigma) f(\sigma)e^{-\vp(\sigma)-\varepsilon|\sigma|^2},
\end{equation}
for all functions $f\in \BF^{\infty,0}_{\vp+\varepsilon|z|^2}$
and such that $\sum |k(z,\sigma)|\le K$ uniformly in $z$. We fix
$p\in[1,\infty)$.  For an arbitrary $g\in \BF^p_\vp$ and $z\in
\C^n$,
$$
f(w)=g(w)e^{2\varepsilon w\cdot \bar z-\varepsilon |z|^2}
$$
belongs to $\BF^{\infty,0}_{\vp+\varepsilon|z|^2}$ and thus we
may apply \eqref{representacio} to obtain
\[
g(z)e^{-\vp(z)}= f(z)e^{-\vp(z)-\varepsilon|z|^2}=
\sum_{\sigma\in\Sigma}k(z,\sigma)
f(\sigma)e^{-\vp(\sigma)-\varepsilon|\sigma|^2}.
\]
Thus
\[
|g(z)|e^{-\vp(z)}\le \sum_{\sigma\in\Sigma}|k(z,\sigma)|
|g(\sigma)|e^{-\vp(\sigma)}e^{-\varepsilon|z-\sigma|^2}.
\]
This together with the inequality $\sum |k(z,\sigma)|\le K$
implies that
\[
\int_{\C^n} |g(z)|^pe^{-p\vp(z)}\lesssim \sum_{\sigma\in\Sigma}
|g(\sigma)|^pe^{-p\vp(\sigma)},
\]
and that
\[
\sup |g(z)|e^{-\vp(z)} \lesssim \sup_{\sigma\in\Sigma}
|g(\sigma)|e^{-\vp(\sigma)}.
\]
\end{proof}
\begin{lem}\label{antidiscret}
Let $W$ be a uniformly flat hypersurface. Let $\Sigma$ be a
uniformly separated sequence contained in $W$. If $\Sigma$ is a
sampling sequence for $\BF^p_\vp$ then $W$ is a sampling
hypersurface for $\BF^p_\vp$.
\end{lem}
\begin{proof}
We only need to prove that for any $z\in W$, the inequality
\begin{equation}\label{mitja}
|f(z)|^pe^{-p\vp(z)}\le C \int_{D_z} |f(x)|^pe^{-p\vp(x)}
\omega^{n-1}(x),
\end{equation}
holds, where $D_z=W\cap B(z,\varepsilon)$, and the constant $C$
may depend on the radius $\varepsilon$ of the ball but not on the
center $z$.  For if \eqref{mitja} holds then for any function
$f\in \BF^p_\vp$,
\[
\begin{split}
\|fe^{-\vp}\|^p_p \lesssim \sum
|f(\sigma)|^pe^{-p\vp(\sigma)}\lesssim \sum_{\sigma\in \Sigma}
\int_{D_\sigma} |f|^pe^{-p\vp} \omega^{n-1} \le \int_W
|f|^pe^{-p\vp} \omega^{n-1}.
\end{split}
\]
In order to prove \eqref{mitja} we need the hypothesis that
$i\ddbar \vp\simeq \omega$. Under this hypothesis we may again invoke 
the existence of a non vanishing function $h\in \so (B(z,\varepsilon))$
such that $e^\vp\simeq |h|$ in $B(z,\varepsilon)$ with constants
independent of $z$. Thus, we may replace $e^{-\vp}$ by $h^{-1}$
in \eqref{mitja} and get the result if we prove that
\[
|g(z)|^p\le C \int_{D_z} |g(x)|^p \omega^{n-1}(x).
\]
If $D_z$ is a hyperplane then the latter estimate holds for all
holomorphic functions $g$ by the sub-mean value property.  In a
general uniformly flat hypersurface the estimate holds because the
distortion introduced in $\omega^{n-1}$ upon rectifying $D_z$ by a
change of variables is uniformly bounded due to property (B) in
Lemma \ref{flat-properties} for uniformly flat hypersurfaces.
\end{proof}

Finally, we are ready to prove Theorem \ref{samp-thm}.   To this
end, let $\varepsilon>0$ be such that $D^-_\vp(W)>1+\varepsilon$.
We start by proving that $W$ is a sampling manifold for
$\BF_{\vp_\varepsilon}^\infty$, where
$\vp_\varepsilon=\vp+\varepsilon |z|^2$. In order to do so, we
use Lemma~\ref{beurlinglem}. We need to check that for any pair
$(V,\psi)\in K^*(W,\vp_\varepsilon)$ the pair $(V,\psi)$ is
determining. This is true in view of
Lemmas~\ref{estabilitatdensitat} and \ref{jensen}. Now we take the
sequence $\Sigma\subset W$ given by Lemma~\ref{discretitzat}. This
sequence $\Sigma$ is a sampling sequence for
$\BF^\infty_{\vp_\varepsilon}$ and thus it is also sampling for
$\BF^p_\vp$ by Lemma~\ref{inclusions}. Finally by
Lemma~\ref{antidiscret} we conclude that $W$ is a sampling
manifold for $\BF^p_\vp$. \qed

\section{An application to sequences in higher dimensions}
\label{appl-section}

Let $\varphi$ be a plurisubharmonic function in $\C ^n$ such that for
some $c > 0$
\[
c \omega \le \ii \di \dbar \varphi \le \frac{1}{c} \omega.
\]
Let $\Gamma$ be a uniformly separated sequence of points in
$\C^n$. We consider the space
\[
\ell ^p _{\varphi} (\Gamma) := \left \{ \{a_{\gamma} \} _{\gamma \in
\Gamma} \subset  \C\ ; \ \sum _{\Gamma} |a_{\gamma} |^p
e^{-p\varphi(\gamma)} < +\infty \right \}.
\]
Recall that $\Gamma$ is an interpolation sequence if for each $\{
a_{\gamma}\} \in \ell ^p _{\varphi}(\Gamma)$ there exists $F \in
\BF_{\varphi}^p (\C ^n)$ such that
\[
F (\gamma ) = a_{\gamma},\quad \gamma \in \Gamma,
\]
and that $\Gamma$ is a sampling sequence if there is a constant $M
>1 $ such that for all $F  \in \BF_{\varphi}^p (\C ^n)$
\begin{eqnarray*}\label{samp-ineq-seq}
\frac{1}{M}\int _{\C ^n} |F|^p e^{-p\varphi} \omega ^{n} \le \sum
_{\Gamma}  |F(\gamma)|^p e^{-p\varphi(\gamma)} \le M \int _{\C ^n}
|F|^p e^{-p\varphi} \omega ^{n}.
\end{eqnarray*}

Sufficient conditions are known for a sequence to be
interpolating, and also sampling.  There are also (different)
necessary conditions.  However, all the known conditions involve
only the number of points of the sequence contained in a large
ball.  It has been known for some time that such a condition could
not possibly characterize interpolation and sampling sequences,
since it does not take into account how points are distributed
relative to one another.  For example, consider the situation of
interpolation. If all the points of a sequence lie on a line, then
to be interpolating there must be at most $O(r^2)$ points in any
ball of radius $r$.  On the other hand, the number of points of a
lattice in $\C ^n$ lying inside a ball of radius $r$ is
$O(r^{2n})$.  Thus any condition for interpolation that takes into
account only the number of points of the sequence lying in a ball
of radius $r$ would not suffice to conclude that any lattice, no
matter how sparse, is an interpolation sequence. Similar reasoning
shows that analogous problems arise in the case of sampling
conditions.

The present paper and the paper \cite{sv} suggest an approach to
studying interpolation and sampling sequences by induction on
dimension.  In \cite{sv} two of us tackled the 1-dimensional case.
The present paper tackles the problem from the other end.  In this
section, we show that the results of the present paper already
improve what is known for sequences in higher dimension.

\subsection{Applications to interpolation}

For simplicity, we restrict to the case of sequences in $\C ^2$.
As mentioned, at present rather poor density conditions are known
in the general higher dimensional case.  However, in a very
symmetric situation there is a characterization of interpolation
and sampling sequences in $\C ^2$.  Suppose the sequence $\Gamma$
is of the form
\[
\Gamma = \Gamma _1 \times \Gamma _2,
\]
where $\Gamma _1, \Gamma _2$ are sequences in $\C$.  Suppose,
moreover, that the weight $\varphi$ splits:
\[
\varphi (z,w) = \varphi _1 (z) + \varphi _2 (w),
\]
where $\Delta \varphi _j \simeq 1$, $j =1,2$.  Then the following is
true:

\begin{claim} $\Gamma$ is interpolating (resp. sampling) with respect
to the weight $\varphi$ if and only if for both $j=1$ and $2$, $\Gamma
_j$ is interpolating (resp. sampling) for the weight $\varphi _j$.
\end{claim}

\noi This result can be recovered from the 1-dimensional
characterization of interpolation and sampling established in \cite{quimbo} and
\cite{quimseep}.

We shall now generalize this result to the case of arbitrary
sequences lying on a family of parallel lines in $\C ^2$.  To this
end, let $\Gamma = \{ \gamma _j\}, \Lambda _1 = \{ \lambda _{1,j}
\},\Lambda _2 = \{ \lambda _{2j}\},\ldots$ be sequences in $\C$.
Define
\[
\Sigma = \left \{ (\gamma _j , \lambda _{jk})\ ;\ j,k=1,2,\ldots
\right \}.
\]
As a corollary of our main results, we have the following theorem.

\begin{thm}\label{seq-interp-thm}
Suppose that for some fixed $\varepsilon > 0$, each $\Lambda _j$ has
density $\le 1- \varepsilon$, with respect to the weight $\varphi (\gamma _j ,
\cdot)$, and that
\begin{eqnarray}\label{split-density}
\frac{\# \Gamma \cap \D (z,r)}{r^2 \Delta _z \varphi (z,w)} <
\frac{\det(\ii \di \dbar \varphi (z,w))}{\Delta _z \varphi
(z,w)\Delta _w \varphi (z,w)}
\end{eqnarray}
for all $z,w \in \C$.  Then $\Sigma$
is interpolating for $\BF_{\varphi}^p (\C ^2).$
\end{thm}

\begin{proof}
Let $W = \Gamma \times \C$.  We first calculate the density of
$W$.  To this end, let $T(z,w)=\sigma (z)$, where $\sigma$ is a
holomorphic function whose zero set, counting multiplicity, is
$\Gamma$.  Then the zero set of $T$ in $\C ^2$ is $W$, and one
sees easily that
\begin{eqnarray*}
D(W,x,r) &=& \sup _{t \in \C} \frac{\sum _j \area \left (
\left ( \{\gamma _j\} \times \C \right )  \cap B(x,r) \right
)}{\vol(B(x,r)) ( \Delta _z \varphi(x) + \Delta _w \varphi (x)
|t|^2 + 2 \Re(\varphi _{z\bar
w} \bar t))}\\
&=& \frac{\sum _j \area \left ( \left ( \{\gamma _j\} \times
\C \right ) \cap B(x,r) \right )}{\vol(B(x,r)) \left (
\Delta _z \varphi(x)
-\frac{|\varphi _{z\bar w}(x)|^2}{\Delta _w \varphi (x)}\right )}\\
&=& \frac{\sum _j \area \left ( \{\gamma _j\} \times \C \cap
B(x,r) \right )}{\vol(B(x,r)) \Delta _z \varphi(x)}\frac{\Delta
_z \varphi(x)\Delta _w\varphi(x)}{\det\left ( \ii \di \dbar \varphi (x)
\right )}.
\end{eqnarray*}
Since we are going to take $\limsup$ as $r \to \infty$, condition
\eqref{split-density} implies that $W$ is an interpolation
hypersurface.

Now suppose given a sequence of values $\{ a_{jk} \}$ such that
\[
\sum _{j} \sum _{k} |a_{jk}|^p e^{p\varphi (\gamma _j,\lambda _{jk})}
<+\infty.
\]
Fix $j$.  Since $\Lambda _j$ is interpolating, there is a function
$g_j (w)$ such that
\[
g_j (\lambda _{jk}) = a_{jk} \quad \text{and} \quad \int _{\C}
|g_j(w)|^p e^{-p\varphi (\gamma _j,w)} dA(w) \le C \sum _k |a_{jk}|^p
e^{p\varphi (\gamma _j,\lambda _{jk})}
\]
for some absolute constant $C$.  (This is not immediate; one has
to use the fact that an interpolation operator can be constructed
with norm depending only on the density of the sequence.  The
uniformity of $C$ now follows because the density of
$\Lambda _j$ is bounded away from 1 uniformly in $j$.)

Define the function $f \in \co (W)$ by
\[
f(\gamma _j ,w) = g_j (w).
\]
Then the estimates on the $L^p$ norms of $g_j$ imply that $f \in
\lf^p _{\varphi} (W).$  By Theorem~\ref{interp-thm}, there exists
$F \in \BF_{\varphi}^p (\C ^2)$ such that $F |W =f$. Thus

\[
F(\gamma _j ,\lambda _{jk}) = f(\gamma _j ,\lambda _{jk}) = \{a_{jk}\},
\]
and the proof is complete.
\end{proof}

We note that, unlike the case of lattices mentioned above, the
condition \eqref{split-density} is not necessary in general, even
for sequences that lie on parallel lines.  To see this, consider
the weight $\varphi (z,w) = |z|^2 + |z+w|^2$.  Let $\Sigma = \{0\}
\times \Gamma$, where $\Gamma$ is a sequence with density between
$\frac{1}{2}$ and $1$.  Then $\Gamma$ is interpolating in $W =
\{0\} \times \C$ and $W$ is interpolating in $\C ^2$.  (In fact,
the density of $W$ is zero.)  But the reader can check that
condition \eqref{split-density} does not hold.  This observation
suggests that perhaps the previously mentioned inductive approach
is lacking another, possibly deep ingredient.

\subsection{Application to sampling sequences}

Let $\Sigma$ be a sequence of the form described before the
statement of Theorem~\ref{seq-interp-thm}.  By analogy with
Theorem~\ref{seq-interp-thm}, we have the following application of
Theorem~\ref{samp-thm} to sequences.

\begin{thm}\label{seq-samp-thm}
Suppose that for some fixed $\varepsilon > 0$, each $\Lambda _j$ has
density $\ge 1+\varepsilon$ with respect to the weight $\varphi
(\gamma _j , \cdot)$ and that, for some $r > 0$,
\begin{eqnarray}\label{split-density-samp}
\frac{\# \Gamma \cap \D (z,r)}{r^2 \Delta _z \varphi (z,w)} >
\frac{\det(\ii \di \dbar \varphi (z,w))}{\Delta _z \varphi
(z,w)\Delta _w \varphi (z,w)}
\end{eqnarray}
for all $z,w \in \C$.  Then $\Sigma$
is sampling for $\BF^p
_{\varphi} (\C ^2).$
\end{thm}

\begin{proof}
Let $W = \Gamma \times \C$.  The upper sampling inequality
holds since $W$ is uniformly flat and $\Sigma \subset W$ is
uniformly separated on each line of $W$.

Next, let $F \in \BF^p_{\varphi} (\C ^2)$. Condition
\eqref{split-density-samp} implies that $W$ is sampling, and thus
\[
\int _{\C ^2} |F |^p e^{-p\varphi} \omega ^2 \le C_1 \int _W
|F|^p e^{-p\varphi}\omega.
\]
Now, since each $\Lambda _j$ is sampling with density bounded away
from $1$ uniformly in $j$, we see that there is $C>0$ such that
for each $j$,
\[
\int _{\{\gamma _j \} \times \C} |F (\gamma _j,w)|^p
e^{-p\varphi(\gamma _j,w)} dA(w) \le C \sum _{k} |F (\gamma
_j,\lambda _{jk})|^p e^{-p\varphi (\gamma _j,\lambda _{jk})}.
\]
Summing over $j$, we have
\[
\int _W |F|^p e^{-p\varphi}\omega \le C_2 \sum _{j,k} |F (\gamma
_j,\lambda _{jk})|^p e^{-p\varphi (\gamma _j,\lambda _{jk})}.
\]
This completes the proof.
\end{proof}

\end{document}